\documentclass[11pt]{article}

\usepackage{setspace}

\usepackage{empheq}
\usepackage{amsmath}
\usepackage{amssymb}
\usepackage{amsthm}
\usepackage{titlesec}
\usepackage{graphicx}
\usepackage{caption}
\usepackage{subcaption}

\usepackage{diagbox}

\usepackage[T1]{fontenc}
\usepackage[utf8]{inputenc}

\usepackage{comment}

\usepackage{indentfirst}
\usepackage[top=1.25in,left=1.25in,right=1.25in]{geometry}
\newlength{\numone}
\newlength{\widone}
\newlength{\numtwo}
\newlength{\widtwo}
\newtheorem{thm}{Theorem}[section]
\newtheorem{lemma}[thm]{Lemma}
\newtheorem{prop}[thm]{Proposition}
\newtheorem{tra}{Transformation}

\numberwithin{equation}{section}

\author{\Large{Riccardo W. Maffucci}}
\newcommand{\Addresses}{{
		\footnotesize
		
		R.W.~Maffucci, \textsc{Dipartimento di Matematica, Universit\`a di Torino\\Via Carlo Alberto 10, Turin 10123, Italy}\par\nopagebreak\vspace{-0.35cm}
		\textit{E-mail address}, R.W.~Maffucci: \texttt{riccardowm@hotmail.com}}}

\title{\Large{\uppercase{\bf On the faces of unigraphic $3$-polytopes}}}

\date{}

\def\calP{\mathcal{P}}

\newcommand{\T}{\mathcal{T}}
\newcommand{\U}{\mathcal{U}}

\begin{document}
\titleformat{\section}
  {\Large\scshape}{\thesection}{1em}{}
\titleformat{\subsection}
  {\large\scshape}{\thesubsection}{1em}{}
\maketitle
\Addresses


\begin{abstract}
A $3$-polytope is a $3$-connected, planar graph. It is called unigraphic if it does not share its vertex degree sequence with any other $3$-polytope, up to graph isomorphism. The classification of unigraphic $3$-polytopes appears to be a difficult problem.

In this paper we prove that, apart from pyramids, all unigraphic $3$-polytopes have no $n$-gonal faces for $n\geq 10$. Our method involves defining several planar graph transformations on a given $3$-polytope containing an $n$-gonal face with $n\geq 10$. The delicate part is to prove that, for every such $3$-polytope, at least one of these transformations both preserves $3$-connectivity, and is not an isomorphism.
\end{abstract}

\section{Introduction}
In this paper, we consider simple graphs with no multiple edges or loops. A graph is the $1$-skeleton of a $3$-polytope (polyhedron) if and only if it is $3$-connected and planar (Rademacher-Steinitz Theorem). We call these {\em polyhedral graphs} or polyhedra. They have been extensively studied by various authors, including e.g.~Gr\"unbaum \cite{grun67}, Klee \cite{klee64}, Barnette \cite{barn73}, Jendrol' \cite{jend97}. A pair of homeomorphic polyhedra corresponds naturally to a pair of isomorphic $3$-connected, planar graphs.

In a polyhedron, the number of vertices $p$, edges $q$, and regions $r$ are related by Euler's formula for planar graphs $p-q+r=2$. The regions are also called faces. They are bounded by polygons, i.e. $n$-cycles with $n\geq 3$ (this fact is true more generally for the regions of any $2$-connected, planar graph).

In fact, a planar graph is $2$-connected if and only if each of its regions is delimited by a cycle. A $2$-connected, planar graph is $3$-connected if, however we choose two distinct regions, their intersection is either empty, or a single vertex, or an edge. 
Two polyhedral faces sharing an edge will be called `adjacent', consistent with the corresponding dual vertices being adjacent in the dual polyhedron \cite[\S 3-4]{dieste}.

For a graph $G$ we may write the degree sequence
\[\sigma: d_1,\dots,d_p\]
in non-increasing order, where $V(G)=\{v_1,\dots,v_p\}$ and $d_i=\deg(v_i)$ for $1\leq i\leq p$. 
If there is only one realisation of a degree sequence up to isomorphism, then the sequence and corresponding graph are called unigraphic. The study of unigraphs has attracted recent attention, e.g. \cite{tysh00,borri9,barr13,hori19,bar023}.

Given a subclass of graphs, one problem is to determine which ones are unigraphic (with respect to the class of all graphs). For the subclass of polyhedra, there are exactly eight solutions. Indeed, more generally, there are exactly eight degree sequences that are {\em forcibly polyhedral}, in the sense that every realisation is polyhedral \cite{mafrao}. Each of these eight sequences has just one realisation.

\begin{figure}[h!]
	\centering
	\includegraphics[width=4.5cm]{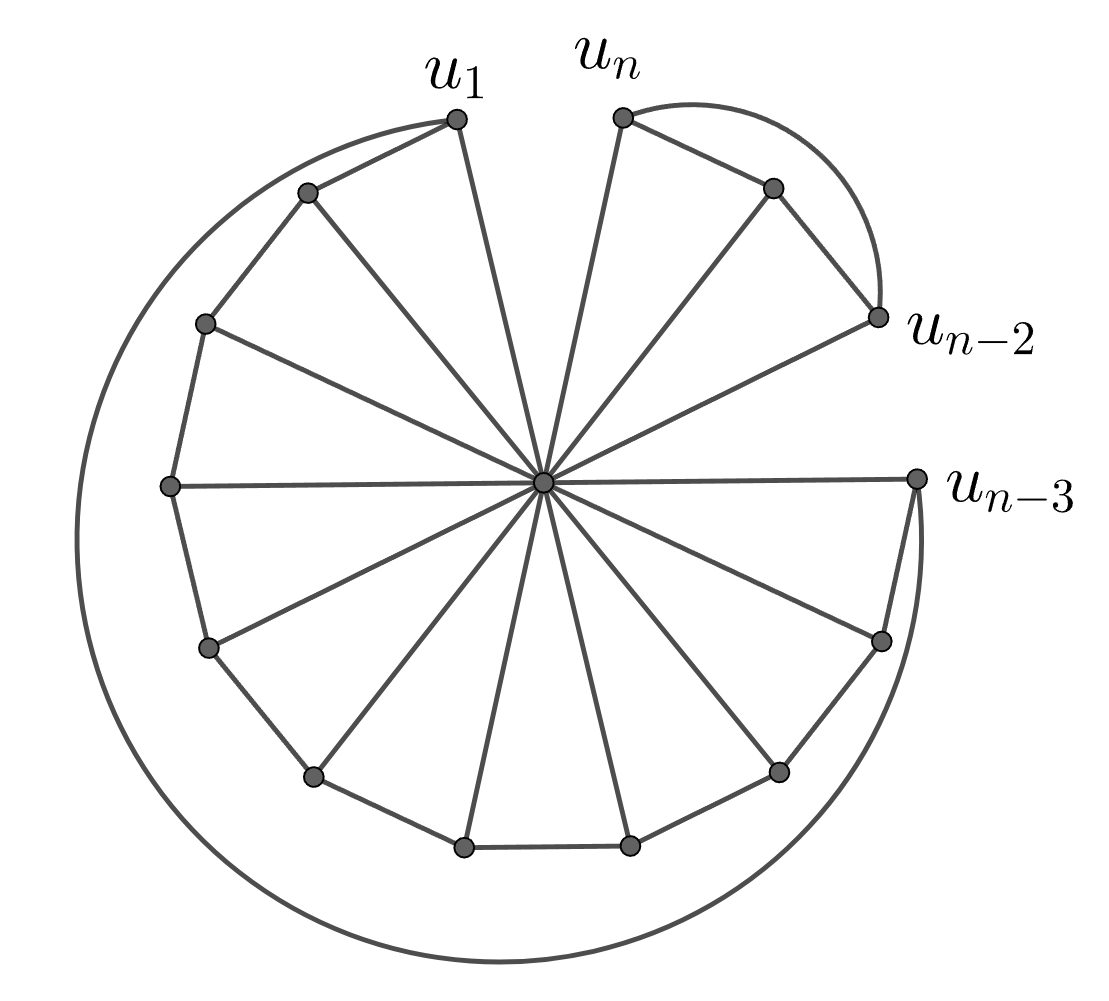}
	\caption{Given the $n$-gonal pyramid $G$ of base $[u_1,u_2,\dots,u_n]$, with  $n\geq 6$, we take $G'=G-u_1u_{n}-u_{n-3}u_{n-2}+u_{1}u_{n-3}+u_{n-2}u_{n}$, obtaining another graph of same sequence as $G$.}
	\label{fig:notpyr}
\end{figure}

Another problem is to characterise the sequences that have a unique realisation {\em within a given subclass of graphs}. For the subclass of polyhedra, this question was posed in \cite{mafp05}. In other words, these are the degree sequences $\sigma$ that have {\bf one and only one polyhedral realisation} (up to graph isomorphism). Henceforth, we will call these sequences, and their corresponding graphs, `polyhedral unigraphic' or `unigraphic' interchangeably. That is to say, henceforth the term `unigraphic' will refer to sequences with exactly one polyhedral realisation. An example is given by the pyramids (note that the $n$-gonal pyramid for $n\geq 6$ is unigraphic w.r.t. the subclass of polyhedra, but not w.r.t. the class of all graphs -- Figure \ref{fig:notpyr}). A few non-trivial examples of families of unigraphic polyhedra were found in \cite{mafp05,delmaf,mafsd1}. In general, this seems to be a difficult problem.

The main theorem of this paper constitutes a major step forward.

\begin{thm}
	\label{thm:1}
	Let $G$ be a unigraphic polyhedron that is not a pyramid. Then $G$ does not contain any $n$-gonal faces for $n\geq 10$.
\end{thm}

Section \ref{sec:2} is entirely dedicated to the proof of Theorem \ref{thm:1}. The proof involves several ideas and graph transformations.

\paragraph{About the proof.} The general strategy to see that a certain polyhedron $G$ is not unigraphic, is to apply a graph transformation that preserves degree sequence, $3$-connectivity, and planarity, and moreover produces a graph not isomorphic to $G$. One attempt is to add diagonals of the polygonal faces, and delete the same number of edges (this preserves planarity), in such a way that the degree of each vertex is invariant, or more generally the vertex degrees are permuted in $G$. The delicate part is to find in what cases such a transformation preserves $3$-connectivity, and when is the resulting graph isomorphic to $G$. To take all scenarios into account, we will need to define several transformations on a given graph, and then prove that at least one of them has the desired properties.

\paragraph{Discussion.}
Our proof of Theorem \ref{thm:1} is constructive, in the following sense. Given a polyhedron containing an $n$-gonal face for $n\geq 10$, our proof also constitutes an algorithm to find another polyhedral realisation of its degree sequence.
\\
This paper constitutes a major step forwards in understanding polyhedral unigraphic sequences, a question posed and partly answered in \cite{mafp05,delmaf,mafsd1}. It follows from Theorem \ref{thm:1} that, apart from pyramids, all faces of all unigraphic polyhedra have nine or fewer edges. The aim of future work is to more precisely describe these unigraphic polyhedra.
\\
In \cite{mafrao} the author revisited Rao's Theorem on forcibly planar degree sequences, to specifically find all forcibly polyhedral sequences. As opposed to this paper, in \cite{mafrao} one needs to find transformations on a polyhedral graph that will preserve degree sequence, but {\em break} planarity or $3$-connectivity.


\paragraph{Notation.} The letter $G$ will denote everywhere a $3$-connected, planar, $(p,q)$ graph (i.e.~with $p$ vertices and $q$ edges), of degree sequence $\sigma$. The dual polyhedron is $G^*$.
\\
The curly letters $\T,\U$ indicate transformations applied to $G$. 
Sketches that illustrate a transformation $\T$ on $G$ will depict a subgraph of $G$, to show how it is modified by $\T$. 
The operations of removing/adding vertices or edges are represented with $-$ and $+$ symbols respectively. The notation $\simeq$ stands for graph isomorphism. 
\\
In a $2$-connected, planar graph, we write
\[F=[u_1,\dots,u_n], \qquad \text{ where } n\geq 3,\]
to indicate that the region $F$ is delimited by an $n$-cycle ($n$-gon) $u_1,\dots,u_n$. If a certain cycle does not necessarily bound a region, it will be written without square brackets.

\section{Proof of Theorem \ref{thm:1}}
\label{sec:2}

\subsection{Graph transformations}
We begin with the preparatory work of defining appropriate graph transformations. In the following, we will write simply $u_i$ in place of $u_{(i \mod n)}$, and set $u_0:=u_n$.
\begin{tra}
	\label{tra:1}
Let $G$ be a polyhedral graph, and
\[F=[u_1,u_2,\dots,u_n], \quad n\geq 6\]
an $n$-gonal face of $G$. We will denote by $F_i$ the face sharing the edge $u_iu_{i+1}$ with $F$, for $1\leq i\leq n$. We fix the vertices $u_i,u_j$, for some
\[
1\leq i,j\leq n, \qquad 3\leq|j-i|\leq n-3,
\]
and perform
\begin{align*}
&\T_1=\T_1(F,i,j): G \to G',\\
&G'=G-u_iu_{i+1}-u_ju_{j+1}+u_iu_{j+1}+u_ju_{i+1}
\end{align*}
(Figure \ref{fig:tra1}).
\begin{figure}[h!]
	\centering
	\includegraphics[width=3.0cm]{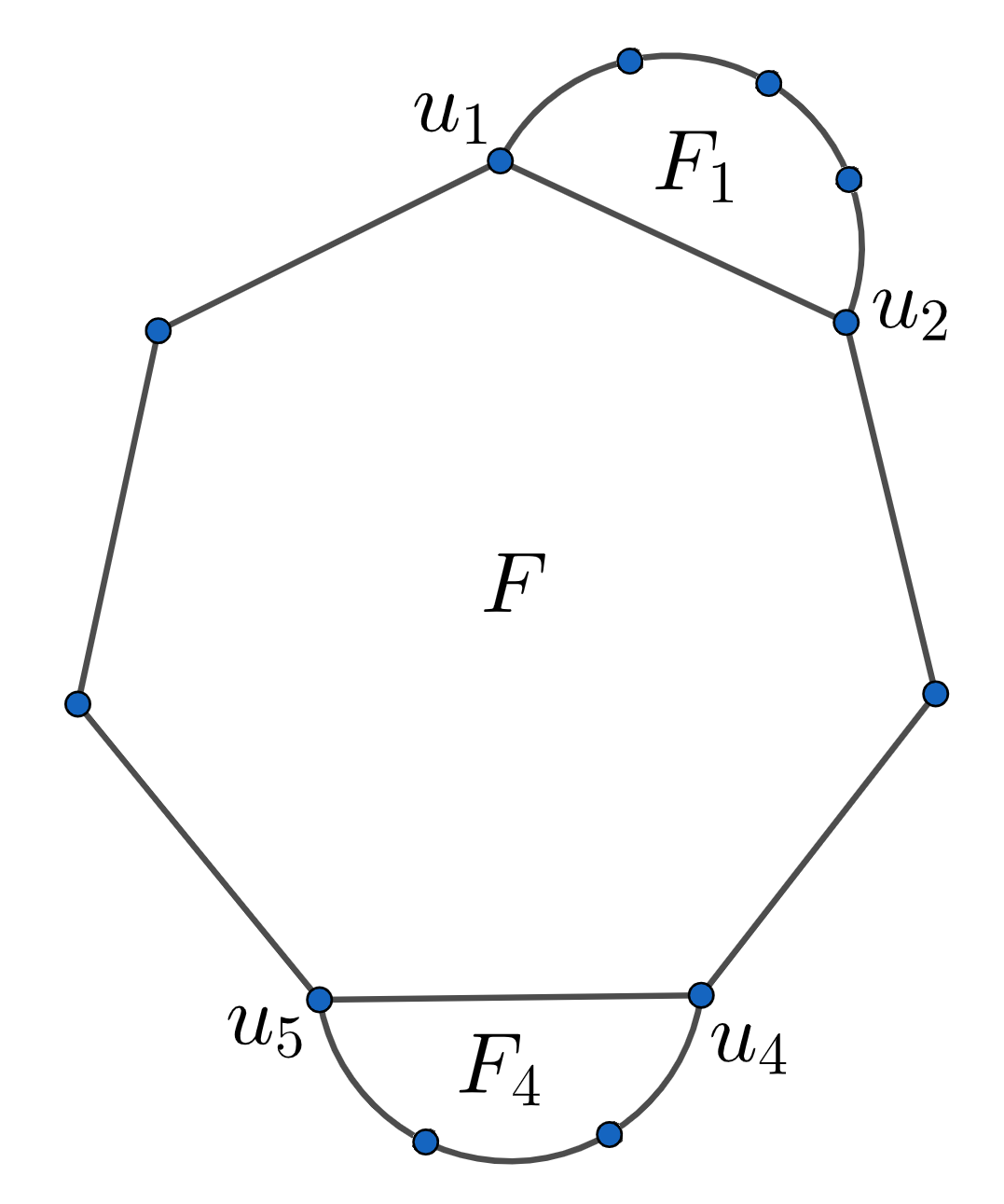}
	\hspace{1.5cm}
	\includegraphics[width=3.0cm]{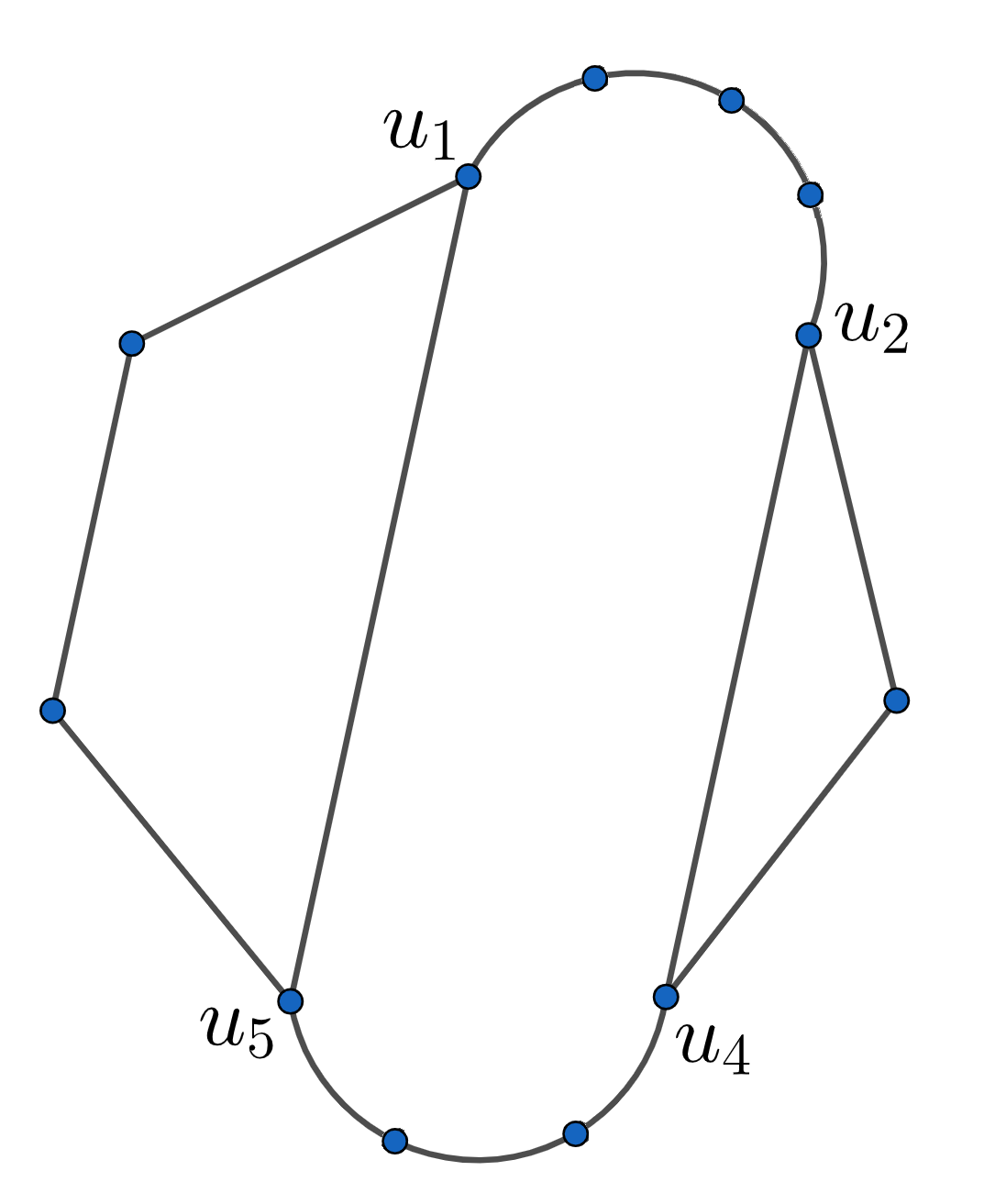}
	\caption{The transformation $\T_1(F,1,4)$ applied to a heptagonal face $F$ (left), transforming the subgraph of $G$ on the left to the one on the right.}
	\label{fig:tra1}
\end{figure}
\end{tra}

Note that in Figure \ref{fig:notpyr}, we have applied $\T_1(F,n-3,n)$ to the pyramid. Let us record a few consideration about Transformation \ref{tra:1}. We start by checking that $G'$ is a simple graph. It suffices to prove that $u_iu_{j+1}\not\in E(G)$. Assume by contradiction that $u_iu_{j+1}\in E(G)$. Since $3\leq|j-i|\leq n-3$, the cycle
\begin{equation}
	\label{eqn:nf1}
u_{j+1},u_{j+2},\dots,u_{i}
\end{equation}
cannot be a face of $G$, otherwise it would share more than two vertices with $F$, contradicting $3$-connectivity. Likewise,
\begin{equation}
\label{eqn:nf2}
u_i,u_{i+1},\dots,u_{j+1}
\end{equation}
is not a face of $G$. We sketch $G$ in the plane. Let $w_1,w_2\in V(G)$ lie inside the cycles \eqref{eqn:nf1} and \eqref{eqn:nf2} respectively. Since $F$ is a face of $G$, every $w_1w_2$-path in $G$ contains at least one of $u_i,u_{j+1}$. That is to say, we have found a $2$-cut $\{u_i,u_{j+1}\}$ in $G$, contradiction.

Clearly $\T_1$ preserves degree sequence and planarity. We also have the following.

\begin{lemma}
	\label{lem:t1}
Let $G,G'$ be as in Transformation \ref{tra:1}. If $G'$ is not $3$-connected, then there exist $x_{i,j},x_{j,i}\in V(G)=V(G')$, such that in $G$, we have $x_{i,j}\in F_i$ and $x_{j,i}\in F_j$, and such that $\{x_{i,j},x_{j,i}\}$ is a $2$-cut in $G'$. In other words, $x_{i,j},x_{j,i}$ lie on the same face ($\neq F$) of $G$. Possibly $x_{i,j}=x_{j,i}$. 
\end{lemma}
\begin{proof}
As mentioned in the Introduction, a planar graph is $2$-connected if and only if each region is delimited by a polygon. Moreover, a $2$-connected, planar graph is $3$-connected if and only if, however we choose two distinct regions, their intersection is either empty, or a vertex, or an edge.

Since Transformation \ref{tra:1} deletes just two edges, the resulting $G'$ is always connected. If $G'$ is of connectivity $1$, then there exists in $G'$ a region $R$ not delimited by a polygon. Since $G$ is $3$-connected, $R$ must be the region of $G'$ containing $u_i,u_{i+1},u_j,u_{j+1}$. Seeing as $G$ is $3$-connected, this can only happen if the faces $F_i,F_j$ of $G$ have non-empty intersection. Therefore, there exists indeed a vertex $x_{i,j}\in F_i\cap F_j$. Note that $x_{i,j}$ lies on every $u_iu_{i+1}$-path in $G'$, hence $x_{i,j}$ is a separating vertex of $G'$ (e.g.~Figure \ref{fig:notpyr}).

Now suppose that $G'$ is of connectivity $2$. This means that there exist two vertices $z_1,z_2\in V(G)$ and distinct regions $R_1,R_2$ in $G'$ such that $V(R_1)\cap V(R_2)\supseteq\{z_1,z_2\}$, with $z_1z_2\not\in E(G')$. Seeing as $G$ is $3$-connected, this can only happen if one of $R_1,R_2$, say $R_1$, is the region of $G'$ containing $u_i,u_{i+1},u_j,u_{j+1}$. Again by $3$-connectivity of $G$, we have $x_{i,j}:=z_1\in F_i$ and $x_{j,i}:=z_2\in F_j$. The region $R_2$ of $G'$ contains both $x_{i,j}$ and $x_{j,i}$. Finally, we note that $R_2$ is also a face of $G$.
\end{proof}

Moreover, if $G'=\T_1(F,i,j)(G)\simeq G$, then either $F_i$ is a $|j-i|$-gon and $F_j$ an $n-|j-i|$-gon, or vice versa (otherwise the new graph would not have the same number of $n$-gons as $G$).

\begin{tra}
	\label{tra:6}
	Let $G$ be a polyhedron, and
	\[F=[u_1,u_2,\dots,u_n],\qquad n\geq 4\]
	an $n$-gonal face of $G$. Assume that the indices $i,j$ satisfy \[d:=\deg(u_i)>\deg(u_j)=:d',\]
	and that
	\[w_1=u_{k},w_2,\dots,w_d=u_{l}\]
	are the vertices adjacent to $u_i$, in order around this vertex in a planar immersion of $G$, with $2\leq|k-j|\leq n-2$. We define
	\begin{align*}
		&\T_2(F,i,k,j):G\to G',\\ &G'=G-u_iw_1-u_iw_2-\dots-u_iw_{d-d'}+u_jw_1+u_jw_2+\dots+u_jw_{d-d'}
	\end{align*}
	(Figure \ref{fig:tra4}).
	\begin{figure}[h!]
		\centering
		\includegraphics[width=4.5cm]{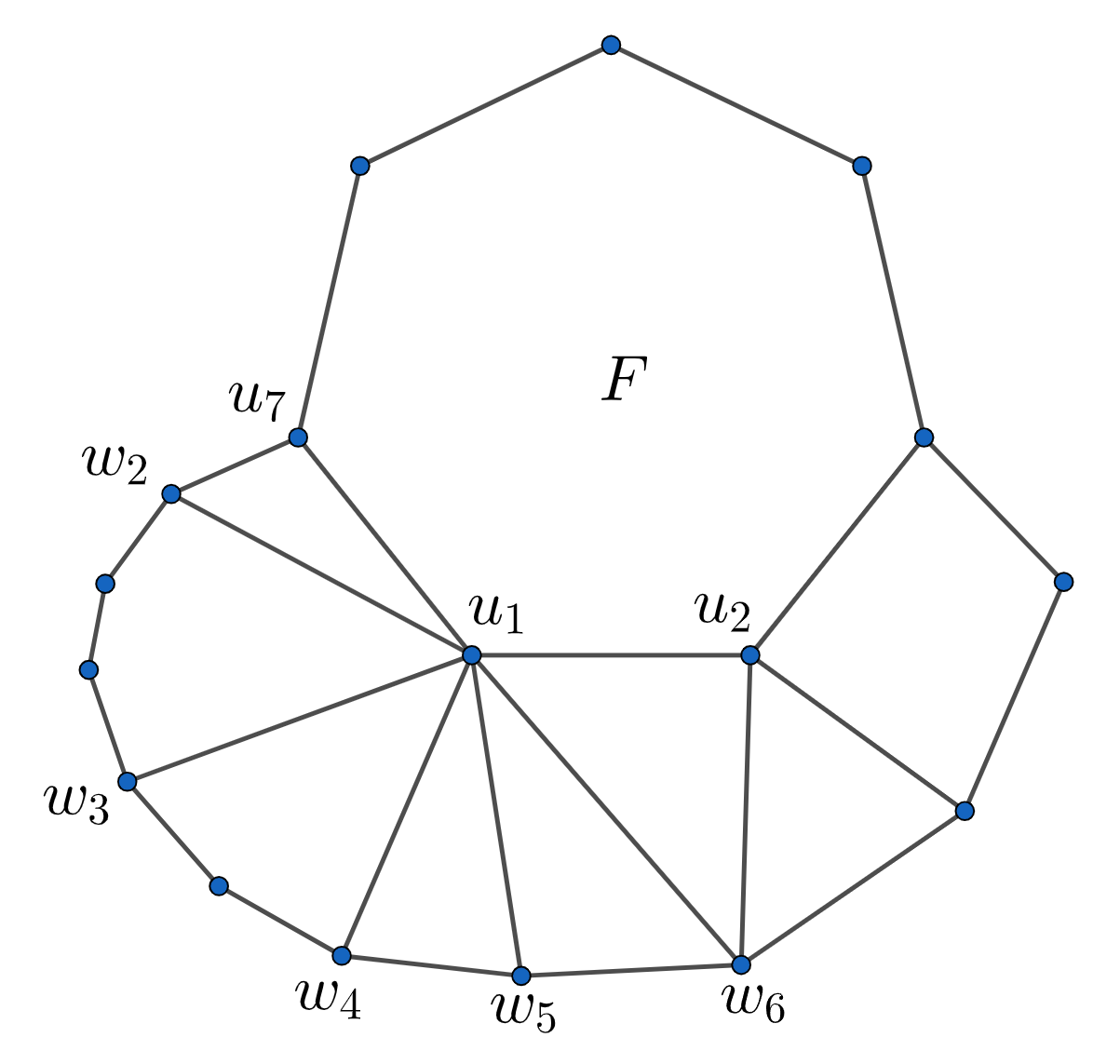}
		\hspace{1.0cm}
		\includegraphics[width=4.5cm]{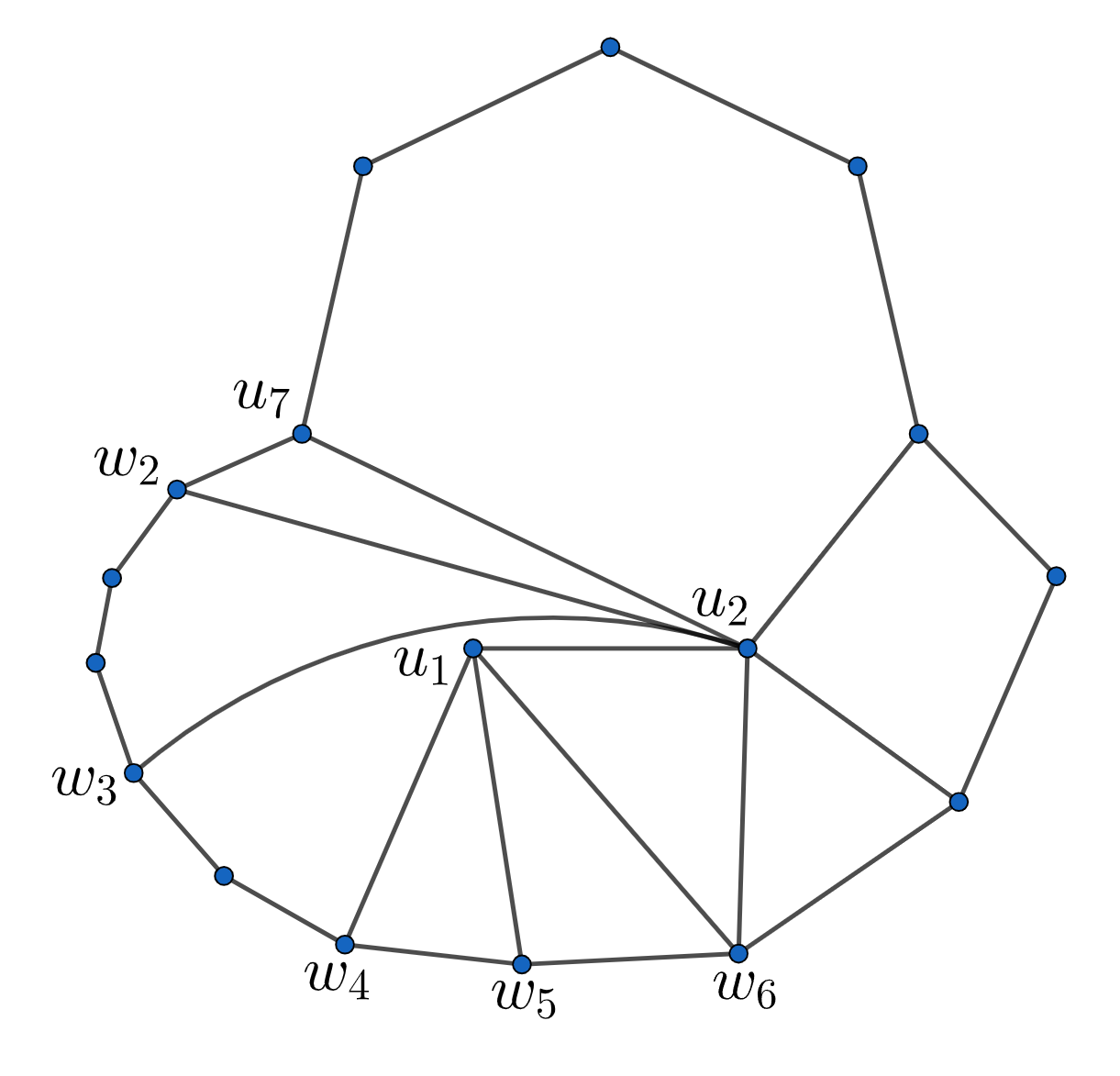}
		\caption{An illustration of the transformation $\T_2(F,1,7,2)$, with $d=7$, $d'=4$, turning the subgraph of $G$ depicted on the left to the one on the right.}
		\label{fig:tra4}
	\end{figure}
	
\end{tra}

Let us make a couple of remarks on Transformation \ref{tra:6}. If $u_j$ is not adjacent in $G$ to any of $w_1,\dots,w_{d-d'}$, then $G'$ is a simple graph. In the proof of Theorem \ref{thm:1}, when performing Transformation \ref{tra:6} on a graph $G$, we will check that $G$ satisfies this condition.

Clearly Transformation \ref{tra:6} preserves planarity and degree sequence. If $G'\simeq G$, then in particular $G,G'$ have the same number of $n$-gonal faces. Thereby, $w_{d-d'},u_i,w_{d-d'+1}$ must be consecutive vertices on the boundary of an $n'$-gonal face in $G$, where
\[n'=\begin{cases}
(j-i)\mod n& \text{ if } k=(i+1)\mod n,\\
n-((j-i)\mod n)& \text{ if } l=(i+1)\mod n.
\end{cases}\]



\begin{tra}
	\label{tra:7}
Let $G$ be a polyhedron, and
\[F=[u_1,u_2,\dots,u_n],\qquad n\geq 4\]
an $n$-gonal face of $G$, with vertices ordered clockwise. We will denote by $F_i$ the face sharing the edge $u_iu_{i+1}$ with $F$, for $1\leq i\leq n$. Fix $i,j$ such that
\[1\leq i,j\leq n, \qquad 2\leq|j-i|\leq n-2.\]
Assume that there exist distinct $x_i\in F_i$, $x_j\in F_j$, such that $x_i,x_j\not\in F$, $x_i,x_j$ lie on the same face, and are not adjacent. Call $a_1$ the clockwise neighbour of $u_i$ on $F_i$, $b_1$ the clockwise neighbour of $x_i$ on the same face, $c_1$ the clockwise neighbour of $u_j$ on $F_j$, and $d_1$ the clockwise neighbour of $x_j$ on the same face. Note that possibly we have one or more of $a_1=x_i$, $b_1=u_{i+1}$, $c_1=x_j$, $d_1=u_{j+1}$. We now transform $G$ as follows,
\begin{align*}
	&\T_3=\T_3(F,i,x_i,j,x_j):G\to G',\\ &G'=G+u_iu_j-u_ia_1+a_1b_1-b_1x_i+x_ix_j-x_jd_1+d_1c_1-c_1u_j
\end{align*}
(Figure \ref{fig:trat}).
\begin{figure}[h!]
	\centering
	\includegraphics[width=3.75cm]{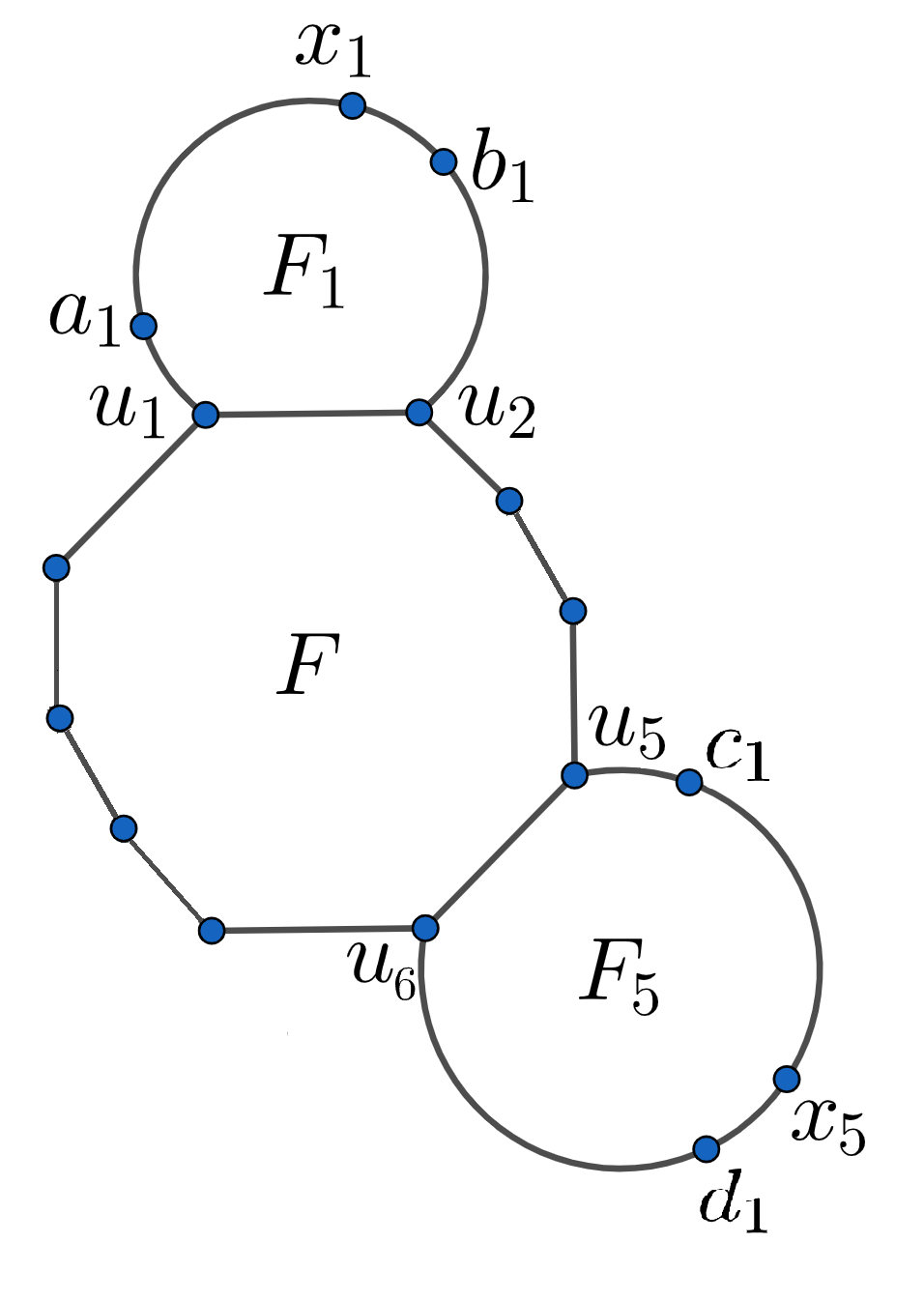}
	\hspace{1.5cm}
	\includegraphics[width=3.75cm]{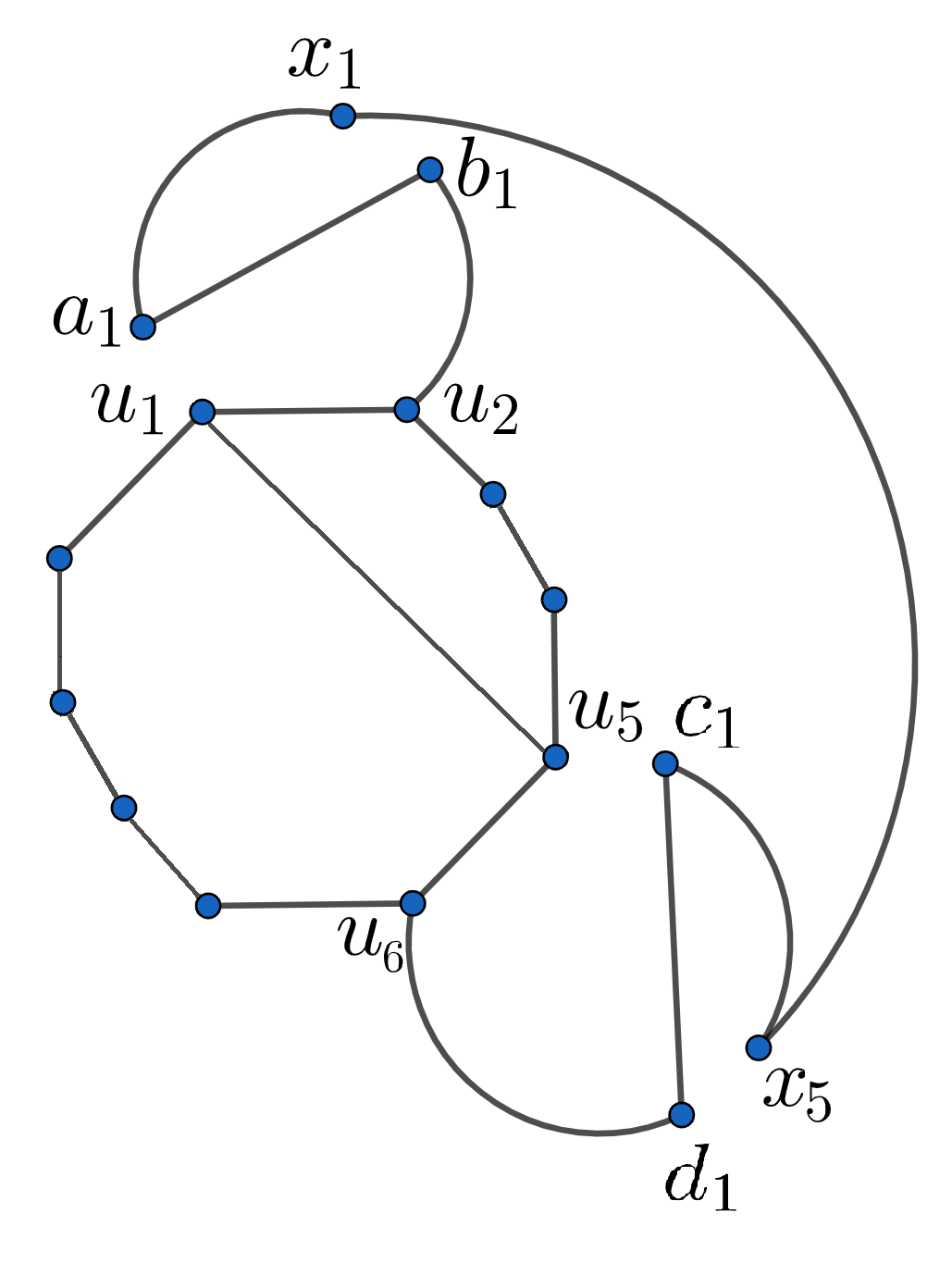}
	\caption{The transformation $\T_3$, in the case of $n=10$, $i=1$, $j=5$, and $F_1,F_5$ bounded by cycles of length at least five. It deletes four edges, and adds the same number.}
	\label{fig:trat}
\end{figure}

\end{tra}
The operation $\T_3$ deletes between two and four edges, adding the same number. It preserves planarity and degree sequence.

\subsection{Setup}
We are ready to begin proving the main theorem. Recall that we write $u_i$ in place of $u_{(i \mod n)}$, where $u_0:=u_n$. 

Let $G$ be a polyhedron of degree sequence $\sigma$, and $F$ an $n$-gonal face
\[[u_1,u_2,\dots,u_n], \quad n\geq 10,\]
where the vertices are ordered clockwise. For $1\leq i\leq n$, we will denote by $F_i$ the face sharing the edge $u_iu_{i+1}$ with $F$.

We consider the transformed graphs
\begin{equation}
	\label{eqn:trgr}
\varphi_1(G),\psi_1(G),\dots,\varphi_n(G),\psi_n(G),
\end{equation}
where
\begin{align}
\label{eqn:zeta}
\varphi_i=\T_1(F,i,(i+\lfloor n/2\rfloor) \mod n),\qquad \psi_i=\T_1(F,i,(i+\lfloor n/2\rfloor-1) \mod n).
\end{align}
These are well-defined instances of $\T_1$, since 
\[\{\lfloor n/2\rfloor,\lfloor n/2\rfloor-1\}\cap\{0,1,2,n-2,n-1\}=\emptyset\]
holds in our case $n\geq 10$. Each graph in \eqref{eqn:trgr} is planar, of same degree sequence as $G$.

If $\varphi_i$ is a graph isomorphism, then either $F_i$ is a $\lfloor n/2\rfloor$-gon and $F_{i+\lfloor n/2\rfloor}$ an $n-\lfloor n/2\rfloor$-gon, or vice versa. Similarly, if $\psi_i$ is an isomorphism, then either $F_i$ is a $\lfloor n/2\rfloor-1$-gon and $F_{i+\lfloor n/2\rfloor-1}$ an $n-\lfloor n/2\rfloor+1$-gon, or vice versa. Therefore, $\varphi_i,\psi_i$ cannot both be graph isomorphisms, otherwise we would get $n+1=2\lfloor n/2\rfloor$, impossible. Now for each $i$, let
\begin{equation*}
\zeta_i=
\begin{cases}
\varphi_i&\varphi_i \text{ is not an isomorphism,}
\\
\psi_i&\text{ otherwise.}
\end{cases}
\end{equation*}
If one of the $\zeta_i$'s preserves $3$-connectivity, then $\zeta_i(G)\not\simeq G$ is a polyhedron with same sequence as $G$, and the proof is complete.

From now on, we thus assume that {\em all} the $\zeta_i$'s break $3$-connectivity. We apply Lemma \ref{lem:t1}. For every $1\leq i\leq n$, since $\zeta_i(G)$ is not $3$-connected, then we can find in $G$ vertices
\begin{equation}
	\label{eqn:x1n}
x_i,y_i
\end{equation}
satisfying $x_i\in F_i$, $y_i\in F_{\xi(i)}$, where
\begin{equation}
	\label{eq:xi}
\xi(i)=\begin{cases}
(i+\lfloor n/2\rfloor)\mod n & \zeta_i=\varphi_i,\\
(i+\lfloor n/2\rfloor-1)\mod n & \zeta_i=\psi_i,
\end{cases}
\end{equation}
and such that $x_i,y_i$ are not both on $F$, and $x_i,y_i$ lie on the same face $F'_i$ in $G$. Note that possibly $x_i=y_i$. To simplify notation, henceforth we will write $x_i$ instead of $x_{(i \mod n)}$, where $x_0:=x_n$, and $y_i$ instead of $y_{(i \mod n)}$, where $y_0:=y_n$.

In the next section, we will prove the following.
\begin{prop}
	\label{prop}
	With the above context and notation, one of the following two statements holds. Either
	\begin{empheq}[box=\fbox]{align}
		\label{eq:ca}
		\notag&\text{there exists a face } F'\neq F \text{ of } G \text{ containing all of } x_1,x_2,\dots,x_n,\\&\text{and moreover, for every } 1\leq i\leq n, \text{ we have } x_i\neq x_{i+2};
	\end{empheq}
	or
	\begin{empheq}[box=\fbox]{align}
		\label{eq:cb}
		\notag&\exists x\in V(G), \ \exists k : 1\leq k\leq n, \quad F_{k+1}=[u_{k+1},u_{k+2},x], 
		\qquad \deg(u_{k+1})=\deg(u_{k+2})=3,
		\\\notag&\text{and moreover,}
		\\&\text{either } \deg(u_k)\geq 4, \qquad \text{ or } F_k \text{ is not a triangle.}
	\end{empheq}
\end{prop}

The rest of the proof of Theorem \ref{thm:1} is organised as follows. In section \ref{sec:prop}, we will prove Proposition \ref{prop}. The two cases \eqref{eq:ca} and \eqref{eq:cb} will be then treated in sections \ref{sec:a} and \ref{sec:b} respectively.

\subsection{Proof of Proposition \ref{prop}}
\label{sec:prop}
We will need a few preliminary lemmas.
\begin{lemma}
	\label{lem:tr}
With the above context and notation, suppose that there exist $1\leq i\leq n$ and $z\in V(G)$ such that
\begin{equation}
z\in F_i,F_{i+1},F_{i+2}.
\end{equation}
Then \eqref{eq:cb} holds.
\end{lemma}
\begin{proof}
Note that $z\neq u_{i+1}$, otherwise the faces $F,F_{i+2}$ would have three common vertices $u_{i+1},u_{i+2},u_{i+3}$. Since the distinct vertices $u_{i+1},z$ belong to both $F_i,F_{i+1}$, we have $u_{i+1}z\in E(G)$. Likewise, the distinct $u_{i+2},z$ belong to both $F_{i+1},F_{i+2}$, hence $u_{i+2}z\in E(G)$. It follows that
\begin{equation}
	\label{eq:j}
\exists j : 1\leq j\leq n, \quad F_{j+1}=[u_{j+1},u_{j+2},z], \quad\text{ and }\quad \deg(u_{j+1})=\deg(u_{j+2})=3,
\end{equation}
namely, $j=i$.

If $\deg(u_{i})\geq 4$ or $F_{i}$ is not a triangle, then we have obtained \eqref{eq:cb} with $k=i$. Otherwise, $\deg(u_{i})=3$ and $F_{i}=[u_{i},u_{i+1},z]$, i.e. \eqref{eq:j} holds with $j=i-1$. If $\deg(u_{i-1})\geq 4$ or $F_{i-1}$ is not a triangle, then we have obtained \eqref{eq:cb} with $k=i-1$. Otherwise, we have \eqref{eq:j} with $j=i-2$. We repeat the same steps until we find $k$ satisfying \eqref{eq:cb}. Such $k$ must exist, since $G$ is not a pyramid.
\end{proof}

\begin{lemma}
	\label{lem:zn}
With the above context and notation, suppose that
\begin{gather}
\label{eq:zam}
\exists z,a,m : z\in V(G), \quad 1\leq a\leq n, \quad 2\leq m\leq n-2,\\\notag z\in F_a,F_{a+m}
\end{gather}
does {\bf not} hold. Then \eqref{eq:ca} holds.
\end{lemma}
\begin{proof}
Recall that, for each $1\leq i\leq n$, we have $x_i\in F_i$, and $y_i\in F_{(i+\lfloor n/2\rfloor)}$ or $F_{(i+\lfloor n/2\rfloor-1)}$. Hence by assumption, the vertices $x_1,x_3,x_5$ are all distinct. Our first goal is to show that
\begin{equation}
	\label{eq:caa}
	\exists k\in\{1,5\} : x_k,x_{3},y_{k},y_{3} \text{ are all distinct, and all lie on the same face } F'\neq F \text{ of } G.
\end{equation}

We have $y_1\in F_{\lfloor n/2\rfloor}$ or $F_{(\lfloor n/2\rfloor+1)}$, and $y_3\in F_{(\lfloor n/2\rfloor+2)}$ or $F_{(\lfloor n/2\rfloor+3)}$. Note that $\lfloor n/2\rfloor\geq 5$, and $\lfloor n/2\rfloor+3\leq n-2$. Therefore, if $y_1\in F_{\lfloor n/2\rfloor}$ or $y_3\in F_{(\lfloor n/2\rfloor+3)}$, then by assumption, since \eqref{eq:zam} does not hold, we obtain \eqref{eq:caa} with $k=1$.

If instead $y_1\in F_{(\lfloor n/2\rfloor+1)}$ and $y_3\in F_{(\lfloor n/2\rfloor+2)}$ hold, then we consider $x_3,x_5,y_3,y_5$. We have $\lfloor n/2\rfloor+2\geq 7$. Moreover, $y_5\in F_{(\lfloor n/2\rfloor+4)}$ or $F_{(\lfloor n/2\rfloor+5)}$, and $\lfloor n/2\rfloor+5\leq n$. Therefore, in this case we obtain \eqref{eq:caa} with $k=5$.

Next, recalling the notation $\xi$ \eqref{eq:xi}, we write
\begin{align*}
	&F_k=[u_k,P_1,
	x_k,\dots,u_{k+1}],
	&F_3=[u_3,P_2,x_{3},\dots,u_{4}
	],
	\\&F_{\xi(k)}=[u_{\xi(k)},P_3,
	y_k,\dots,u_{\xi(k)+1}],
	&F_{\xi(3)}=[u_{\xi(3)},P_4,y_3,\dots,
	u_{\xi(3)+1}],
\end{align*}
where $P_1,P_2,P_3,P_4$ are simple paths. These four faces are pairwise disjoint by assumption.

Recall that there exists a face $F'_k\neq F$ containing $x_k,y_k$, and there exists a face $F'_3\neq F$ containing $x_3,y_3$. We write
\[F'_k=[x_k,P_5,y_k,P_6] \quad \text{ and } \quad F'_3=[x_3,P_7,y_3,P_8],\]
where $P_5,P_6,P_7,P_8$ are simple paths, with $P_5,P_6$ internally disjoint, and $P_7,P_8$ internally disjoint.

Therefore, there exist in $G$ two paths
\begin{equation}
	\label{eqn:paths}
	u_k,P_1,x_k,P_5,y_k,P_3,u_{\xi(k)} \quad \text{ and } \quad u_{3},P_2,x_3,P_7,y_3,P_4,u_{\xi(3)},
\end{equation}
both external to the face $F$. By planarity, the paths \eqref{eqn:paths} must intersect at a vertex $v_1$. Since the faces $F_k,F_3,F_{\xi(k)},F_{\xi(3)}$ are pairwise disjoint, it follows that $P_5$ and $P_7$ intersect at $v_1$. Likewise, $P_5,P_8$ intersect at a vertex $v_2$, $P_6,P_7$ intersect at a vertex $v_3$, and $P_6,P_8$ intersect at a vertex $v_4$. Recall that $P_5,P_6$ are internally disjoint, and $P_7,P_8$ are internally disjoint. Therefore, $F'_k\cap F'_3$ contains four distinct vertices $v_1,v_2,v_3,v_4$, however two distinct faces may intersect in at most two vertices, thus in fact $F'_k=F'_3=:F'$, so that there exists a face $F'\neq F$ of $G$ containing all of $x_k,x_3,y_k,y_3$.

It remains to prove that $F'$ contains all of $x_1,x_2,\dots,x_n$. Recall that the indices are taken modulo $n$, with $x_0=x_n$. We begin by showing that $x_5\in F'$ also when $k=1$. Recall that there exists a face $F'_5\neq F$ containing $x_5,y_5$. In a planar immersion of $G$, the vertices $x_5,y_5$ are cut off from each other by the faces $F,F',F_3,F_{\xi(3)}$. We note that
\[3+2\leq 5\leq \xi(3)-2,\]
so that by assumption $x_5\not\in F_3,F_{\xi(3)}$. The only possibility is thus $x_5,y_5\in F'$, i.e. $F'_5=F'$. Next, $x_7,y_7$ are cut off from each other by $F,F',F_5,F_{\xi(5)}$, and moreover $5+2\leq 7\leq \xi(5)-2$, so that $x_7,y_7\in F'$, and we continue in this fashion proving that all odd labelled vertices $x_1,x_2,\dots,x_n$ belong to $F'$. Now if $n$ is even, we also consider $x_4,y_4$, that are cut off from each other by $F,F',F_3,F_5$. We cannot have $x_4,y_4\in F$, or $x_4,y_4\in F_3$, or $x_4,y_4\in F_5$, so that the only possibility is $x_4,y_4\in F'$. We continue in this way, proving that indeed $F'$ contains all of $x_1,x_2,\dots,x_n$ in this order.
\end{proof}

\begin{lemma}
	\label{lem:zy}
With the above context and notation, assume that
$u_a,u_i,u_b,u_j,u_c,u_k$
appear around the boundary of $F$ in this order, where $1\leq i<j<k\leq n$,
\[j-i\geq 3, \ \text{ and } \ k-j\geq 3.\]
If
\[\exists z\in V(G) : z\in F_a,F_b,F_c, \quad z\not\in F_i,F_j,F_k,\]
then
\[k-j=j-i \ \text{ and } \ n-(k-i)\leq 2.\]
\end{lemma}
\begin{proof}
As remarked after Transformation \ref{tra:1},
\begin{equation}
	\label{eq:zy1}
	\T_1(F,i,j)(G)
\end{equation}
is a simple (since $j-i\geq 3$), planar graph, of same degree sequence $\sigma$ as $G$. Moreover, since $z\in F_a,F_{b}$, $z\not\in F_{i},F_{j}$, and $G$ is planar, then there does not exist a face $\neq F$ of $G$ intersecting both $F_{i},F_{j}$. By Lemma \ref{lem:t1}, \eqref{eq:zy1} is $3$-connected. Thereby, \eqref{eq:zy1} is a $3$-polytope of same sequence as $G$, so that by unigraphicity $\T_1(F,i,j)(G)\simeq G$. In particular, as remarked after Transformation \ref{tra:1}, the face $F_{i}$ is a $j-i$-gon and the face $F_{j}$ is an $n-(j-i)$-gon, or vice versa.

Similarly, $\T_1(F,j,k)(G)\simeq G$, so that $F_{j}$ is a $k-j$-gon and $F_{k}$ an $n-(k-j)$-gon, or vice versa. Hence either $n-(k-j)=j-i$, i.e. $k-i=n$ contradiction, or $k-j=j-i$.

As for the second statement, we argue by contradiction. If $n-(k-i)\geq 3$, then as above $\T_1(F,i,k)$ produces a simple graph, and $\T_1(F,i,k)(G)\simeq G$. It follows that $F_{i}$ is a $k-i$-gon and $F_{k}$ an $n-(k-i)$-gon, or vice versa. On the other hand, as seen above, $F_j$ is a $j-i$-gon, so that $F_i$ is an $n-(j-i)$-gon and $F_k$ is an $n-(k-j)$-gon. It follows that either $n-(k-i)=n-(j-i)$, i.e. $j=k$ contradiction, or $n-(k-i)=n-(k-j)$, i.e. $i=j$ contradiction.
\end{proof}

\begin{proof}[Proof of Proposition \ref{prop}]
Let us assume by contradiction that \eqref{eq:ca} and \eqref{eq:cb} are both false. Then \eqref{eq:zam} holds,
\begin{gather*}
	\exists z,a,m : z\in V(G), \quad 1\leq a\leq n, \quad 2\leq m\leq n-2,\\z\in F_a,F_{a+m}
\end{gather*}
otherwise we would obtain \eqref{eq:ca} immediately via Lemma \ref{lem:zn}. We now ask ourselves how many, and which of the faces
\begin{equation}
	\label{eq:fac}
F_1,F_2,\dots,F_n
\end{equation}
contain $z$. We may always assume that three consecutive elements of \eqref{eq:fac} (regarding $F_n,F_1$ as consecutive) do not all contain $z$, otherwise we would obtain \eqref{eq:cb} via Lemma \ref{lem:tr}.


We claim that at least three pairwise non-consecutive elements of \eqref{eq:fac} contain $z$. We argue by contradiction. Since three consecutive elements of \eqref{eq:fac} cannot all contain $z$, and moreover there do not exist three pairwise non-consecutive elements of \eqref{eq:fac} containing $z$, then in fact at most four elements of \eqref{eq:fac} may contain $z$, and given any three of them, exactly two of them are consecutive. W.l.o.g., $z\not\in F_3,F_4,F_n$, and $z\in F_2,F_{b}$ with $5\leq b\leq n-1$. By the first part of the proof of Lemma \ref{lem:zy}, $\T_1(F,3,n)(G)\simeq G$, so that one of $F_3$, $F_n$ is a triangle and the other an $n-3$-gon, and similarly $\T_1(F,4,n)(G)\simeq G$, so that one of $F_4$, $F_n$ is a quadrangle and the other an $n-4$-gon. This leads to either $n-3=n-4$ or $n-3=4$, contradiction.

Therefore, there exist distinct $i,j,k\in\{1,2,\dots,n\}$ such that
\[z\in F_{i-1},F_{j-1},F_{k-1}, \quad z\not\in F_{i},F_{j},F_{k}.\]
As seen in Lemma \ref{lem:zy}, at least one of
\begin{equation}
	\label{eq:ijk}
|i-j|, \ |j-k|, \ |i-k|
\end{equation}
equals $2$ or less, thus it is exactly $2$.

We consider two cases: either exactly two of \eqref{eq:ijk} equal $2$, or exactly one of them equals $2$. In the first case, we may then write w.l.o.g.
\[z\in F_{n-1},F_{1},F_{3}, \quad z\not\in F_{n},F_{2},F_{4}.\]
Assuming $z\not\in F_5$, we apply Lemma \ref{lem:zy} with $a=1$, $i=2$, $b=3$, $j=5$, $c=n-1$, and $k=n$, to obtain $k-j=j-i$ thus $n=8$, contradiction. Hence $z\in F_5$.

Likewise, assuming $z\not\in F_7$, we apply Lemma \ref{lem:zy} with $a=1$, $i=2$, $b=5$, $j=7$, $c=n-1$, and $k=n$, to obtain $n-7=7-2$. We also apply it with $a=3$, $i=4$, $b=5$, $j=7$, $c=n-1$, and $k=n$, to obtain $n-7=7-4$, contradicting $n-7=7-2$. Hence $z\in F_7$.

At least one of the three consecutive $F_5,F_6,F_7$ does not contain $z$, hence $z\not\in F_6$. We apply Lemma \ref{lem:zy} with $a=1$, $i=2$, $b=5$, $j=6$, $c=n-1$, and $k=n$, obtaining $n=10$. The three consecutive $F_7,F_{8},F_9$ cannot all contain $z$, thus $z\not\in F_{8}$. Summarising,
\[
	n=10, \quad z\in F_1,F_3,F_5,F_7,F_9, \quad z\not\in F_2,F_4,F_6,F_8,F_{10}.
\]
Looking at the isomorphism $\T_1(F,2,6)$, we conclude that $F_2,F_6$ are a $4$-gon and a $6$-gon in some order. If $F_2$ is a $4$-gon and $F_6$ a $6$-gon (the other case is analogous), we analyse in turn the isomorphisms $\T_1(F,6,10)$, $\T_1(F,10,4)$, $\T_1(F,4,8)$, and $\T_1(F,8,2)$, to see that $F_{10}$ is a $4$-gon, $F_{4}$ a $6$-gon, $F_8$ a $4$-gon, and $F_2$ a $6$-gon, contradiction.

It remains to inspect the second case, where exactly one of \eqref{eq:ijk} equals $2$. We write w.l.o.g.
\[z\in F_{n-1},F_{1},F_{m-1}, \quad z\not\in F_{n},F_{2},F_{m},\]
where $5\leq m\leq n-3$ otherwise we would be in the already analysed first case. Applying Lemma \ref{lem:zy} with $a=1$, $i=2$, $b=m-1$, $j=m$, $c=n-1$, and $k=n$, we get $n-m=m-2$ thus $m=n/2+1$, and in particular $m\geq 6$.

Assuming $z\not\in F_3$, we apply Lemma \ref{lem:zy} with $a=1$, $i=3$, $b=n/2$, $j=n/2+1$, $c=n-1$, and $k=n$, to obtain $n-(n/2+1)=(n/2+1)-3$, contradiction. Hence $z\in F_3$.

Next, assuming $z\not\in F_5$, we apply Lemma \ref{lem:zy} with $a=1$, $i=2$, $b=3$, $j=5$, $c=n-1$, and $k=n$, to obtain $n=8$, contradiction. Hence $z\in F_5$.

At least one of the three consecutive $F_3,F_4,F_5$ does not contain $z$, hence $z\not\in F_4$. Now if $n=10$ we conclude as in the previous case. Otherwise, we apply Lemma \ref{lem:zy} with $a=3$, $i=4$, $b=n/2$, $j=n/2+1$, $c=n-1$, and $k=n$, to obtain $n-(n/2+1)=(n/2+1)-4$, contradiction. The proof of Proposition \ref{prop} is complete.
\end{proof}

\subsection{The case \eqref{eq:ca}}
\label{sec:a}
Back to the proof of Theorem \ref{thm:1}, assume for the moment that we are in the case \eqref{eq:ca}, i.e.~there exists a face $F'$ containing all of $x_1,x_2,\dots,x_n$ in this order, and for each $1\leq i\leq n$, we have $x_i\neq x_{i+2}$. We claim that if $F$ contains some of $x_1,x_2,\dots,x_n$, then it contains at most three of them, and these have consecutive indices. Indeed, by contradiction, assume that $x_i,x_{i+3}\in F$. Then we have $x_i\in \{u_i,u_{i+1}\}$ and $x_{i+3}\in \{u_{i+3},u_{i+4}\}$, so that the faces $F,F'$ share two non-adjacent vertices, contradiction.

In particular, w.l.o.g. we may assume the following:
\begin{equation*}
	x_1,x_3,x_5,x_7 \text{ from }\eqref{eqn:x1n} \text{ are all distinct, and moreover } x_1,x_5\not\in F.
\end{equation*}
By planarity, $x_1x_5\not\in E(G)$. We may write in clockwise order
\[F_1=[u_1,a_1,a_2,\dots,a_{\alpha-1},x_1=a_\alpha,b_1,b_2,\dots,b_{\beta-1},u_2=b_\beta]\]
and
\[F_5=[u_5,c_1,c_2,\dots,c_{\gamma-1},x_5=c_\gamma,d_1,d_2,\dots,d_{\delta-1},u_6=d_\delta],\]
for some $\alpha,\beta,\gamma,\delta\geq 1$. Note that possibly $a_1=x_1$, and possibly $c_1=x_5$. We transform $G$ as follows,
\begin{align}
	\notag&\T_3=\T_3(F,1,x_1,5,x_5):G\to G',\\ &G'=G+u_1u_5-u_1a_1+a_1b_1-b_1x_1+x_1x_5-x_5d_1+d_1c_1-c_1u_5.
	\label{eqn:T3}
\end{align}

The resulting graph $G'$ is still planar, and shares its degree sequence with $G$. Let us show $3$-connectivity. We begin by defining in $G$ the faces $A,B$ adjacent to $F_1$, and $C,D$ adjacent to $F_5$, by writing
\begin{align*}
	&A=[u_1,a_1,a'_3,\dots,a'_{l_A}],
	&B=[x_1,b_1,b'_3,\dots,b'_{l_B}],\\
	&C=[u_5,c_1,c'_3,\dots,c'_{l_C}],
	&D=[x_5,d_1,d'_3,\dots,d'_{l_D}].
\end{align*}
Then the regions of $G'$ include
\begin{align*}
	&A'=[a_1,a'_3,\dots,a'_{l_A},u_1,b_{\beta},b_{\beta-1},\dots,b_1],\\
	&B'=[b_1,b'_3,\dots,b'_{l_B},x_1,a_{\alpha-1},a_{\alpha-2},\dots,a_1],
	\\
	&C'=[c_1,c'_3,\dots,c'_{l_C},u_5,d_{\delta},d_{\delta-1},\dots,d_1],\\
	&D'=[d_1,d'_3,\dots,d'_{l_D},x_5,c_{\gamma-1},c_{\gamma-2},\dots,c_1].
\end{align*}
Note that if $a_1=x_1$ (resp. $c_1=x_5$), then $B'=B$ (resp. $D'=D$). By planarity, if $A,B$ (resp. $C,D$) share a vertex, then this vertex is $a_1=x_1$ (resp. $c_1=x_5$). We are ready to show $3$-connectivity of $G'$. First, to see that $G'$ is $2$-connected, one checks that each new region is delimited by a cycle. For instance in $A'$, the relation $a_{\lambda}'=b_{\mu}$ for some $\lambda,\mu$ is not possible, seeing as $G$ is a polyhedron. Similar considerations can be made for $B',C',D'$. For $3$-connectivity of $G'$ to hold, it remains to check that, for each $X'=A',B',C',D'$, the only region of $G'$ containing two non-adjacent vertices of $X'$ is $X'$ itself. By planarity, in $G'$ the only region containing one of
\[a'_3,\dots,a'_{l_A}\]
and also containing one of
\[b_1,b_2,\dots,b_\beta\]
is $A'$ itself. Likewise, the only region of $G'$ that may contain one of
\[b'_3,\dots,b'_{l_B}\]
and also one of the
\[a_1,a_2,\dots,a_{\alpha-1}\]
is $B'$. 
Analogous considerations can be made for $C',D'$. This proves that $G'$ is indeed $3$-connected.

To finish the proof in the case \eqref{eq:ca}, we will now prove that $G'\not\simeq G$. In $G$, the face $F$ satisfies the following property $\calP$: there exists a face $F'$ having non-empty intersection with every face adjacent to $F$. Now $F$ is not a face in $G'$. By contradiction, if $G'\simeq G$, then there must exist a face of $G'$ that is not a face of $G$ and that verifies $\calP$. The only candidates are $A',B',C',D'$, and the two faces created by adding the edge $x_1x_5$. However, by planarity, none of them can have the property $\calP$. The proof of Theorem \ref{thm:1} is complete in this case.

\subsection{The case \eqref{eq:cb}}
\label{sec:b}
Now assume instead that \eqref{eq:cb} holds, so that fixing w.l.o.g. $k=1$ 
we have
\begin{align*}
	\notag&\exists x\in V(G) : \quad F_2=[u_2,u_3,x], \qquad\deg(u_2)=\deg(u_3)=3,
	\\\notag&\text{and moreover,}
	\\&\text{either } \deg(u_1)\geq 4, \qquad \text{ or } F_1 \text{ is not a triangle.}
\end{align*}
If $\deg(u_1)\geq 4$, we consider the transformation
\begin{equation*}
	\T_2(F,1,n,2) : G\to G(2).
\end{equation*}
The vertices adjacent to $u_1$ in $G$ are
\[w_1=u_n,w_2,\dots,w_{\deg_G(u_1)}=u_2,\]
say, ordered clockwise around $u_1$. 
Note that $u_1,x,w_{\deg_G(u_1)-1}$ lie on the same face, and possibly $w_{\deg_G(u_1)-1}=x$. None of 
\[w_1,w_2,\dots,w_{\deg_G(u_1)-3}\]
may coincide with any neighbour of $u_2$ (the neighbours of $u_2$ are $u_1,u_3,x$), hence $G(2)$ is a simple graph. It is planar and shares its sequence with $G$. Let's show that it is also $3$-connected. By contradiction, there exists in $G$ a region containing both $u_2$ and a vertex $c\not\in F$, where $c$ is such that $c,u_1,w_j$ lie on a face in $G$, for some $1\leq j\leq\deg_G(u_1)-3$. This contradicts planarity.

It remains to see what happens if $G(2)\simeq G$. As recorded after defining Transformation \ref{tra:6}, in this case  
in $G$ the vertices
\[u_1,w_{\deg_G(u_1)-3},w_{\deg_G(u_1)-2}\]
lie on an $n-1$-gonal face $F''$, say (otherwise $G$ would have strictly more $n$-gonal faces than the transformed polyhedron). Here we also consider
\begin{equation}
	\label{eqn:T43}
	\T_2(F,1,n,3):G\to G(3).
\end{equation}
As above, $G(3)$ is a simple, planar graph of sequence $\sigma$. It is $3$-connected, otherwise we could find in $G$ a region containing vertices $c',c$, such that $c'\in\{u_2,u_3\}$, while $c,u_1,w_j$ lie on the same face in $G$, for some $1\leq j\leq\deg_G(u_1)-3$, and moreover $c\not\in F$, again impossible. If $G(3)\simeq G$, then $F''$ is an $n-2$-gon. Thereby, at least one of $G(2),G(3)$ is not isomorphic to $G$.

\begin{figure}[h!]
	\centering
	\includegraphics[width=3.75cm]{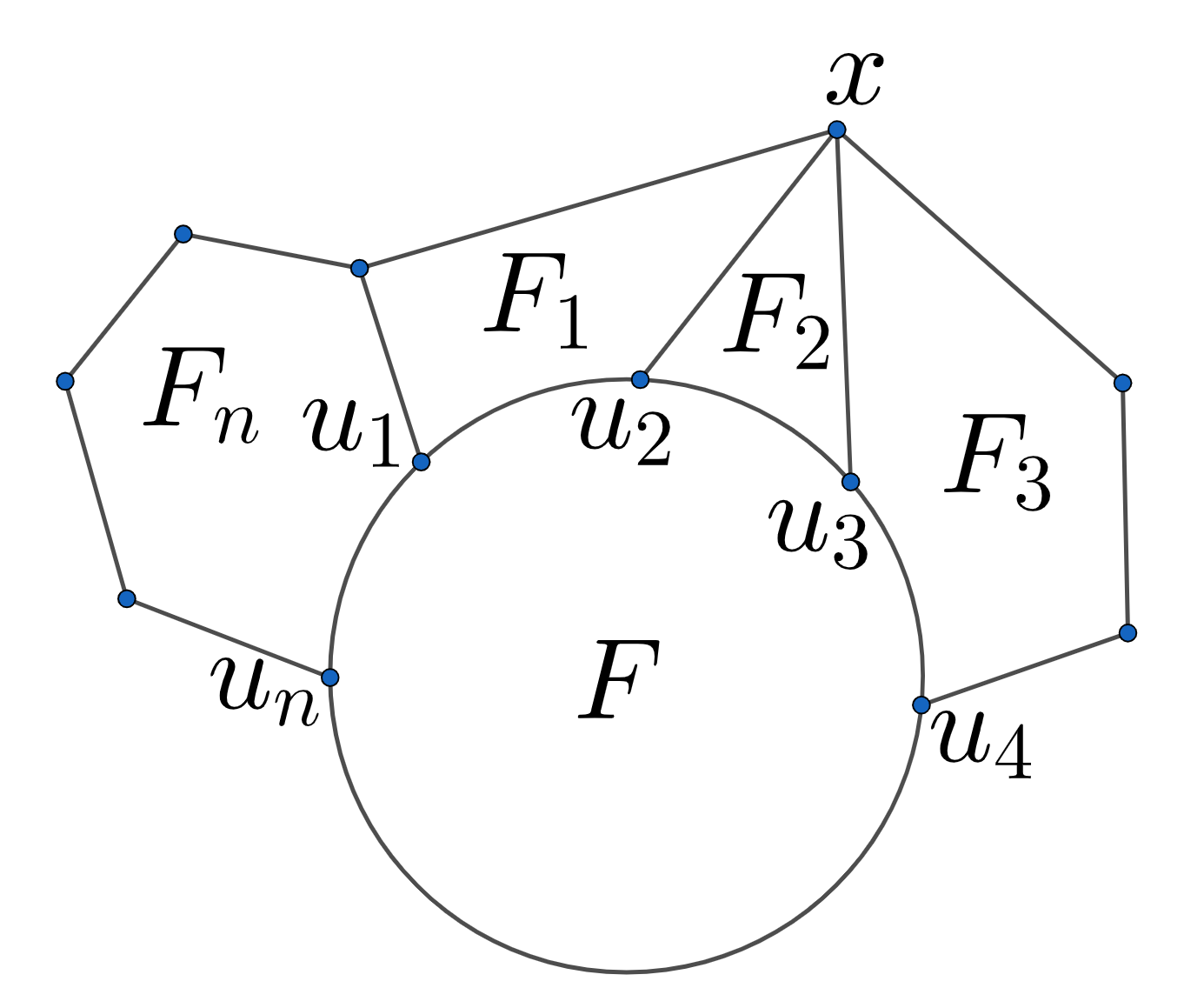}
	\hspace{1.5cm}
	\includegraphics[width=3.75cm]{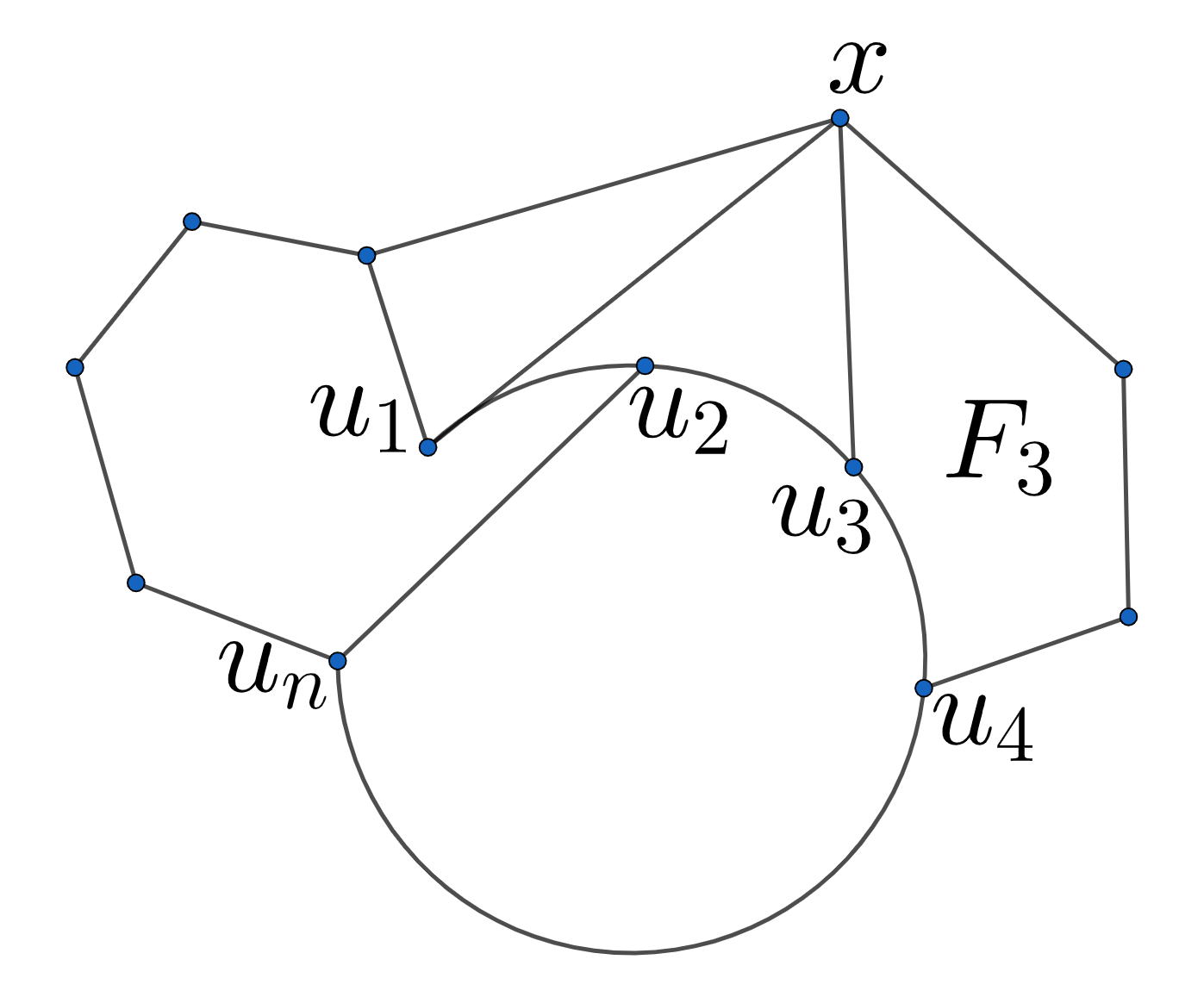}
	\caption{An example of $\U(n,1,2)$, transforming the subgraph of $G$ on the left to the subgraph of $G[n,1,2]$ on the right.}
	\label{fig:G2}
\end{figure}

We now analyse what happens when \eqref{eq:cb} holds, with $k=1$ say, together with
\[F_1=[x,u_2,u_1,\dots]\]
not being a triangle (the dots replace at least one and of course finitely many vertices). Here we consider the transformation
\begin{align*}
	&\U(n,1,2):G\to G[n,1,2],
	\\
	&G[n,1,2]=G+u_1x-xu_2+u_2u_n-u_nu_1,
\end{align*}
illustrated in Figure \ref{fig:G2}. The transformed graph $G[n,1,2]$ is planar and of sequence $\sigma$.

We now show $3$-connectivity for $G[n,1,2]$. It suffices to check that however we choose two distinct regions of $G[n,1,2]$, their intersection is either empty, or a vertex, or an edge. Arguing by contradiction, let $R_1,R_2$ be two regions of $G[n,1,2]$ not satisfying the condition. That is to say, the intersection of $R_1,R_2$ contains two vertices that are not adjacent in  $G[n,1,2]$. Since $G$ is $3$-connected, at least one of $R_1,R_2$, say $R_1$, is not a face of $G$. Hence there are the following two cases. First case, $R_1=[u_1,u_2,u_3,x]$ and $R_2\neq R_1$ contains $u_1$ and $u_3$. However, by \eqref{eq:cb}, the only regions $\neq R_1$ of $G[n,1,2]$ containing $u_3$ are $F_3$ and $[u_2,u_3,\dots,u_n]$, and neither of them contains $u_1$. Second case, $R_1$ is formed by adding $u_2$ to the vertices of $F_n$, and $R_2\neq R_1$ contains $u_2$ and a vertex of $F_n$ distinct from $u_1,u_n$. By \eqref{eq:cb}, the only regions $\neq R_1$ of $G[n,1,2]$ containing $u_2$ are $[u_1,u_2,u_3,x]$ and $[u_2,u_3,\dots,u_n]$, and neither may contain a vertex of $F_n$ distinct from $u_1,u_n$, unless this vertex is $x$. However, in the polyhedron $G$, we have $x\not\in F_n$, otherwise $F_1,F_n$ would share both $u_1,x$, but $u_1x\not\in E(G)$ since $F_1$ is not a triangle. We conclude that $G[n,1,2]$ is $3$-connected. Therefore, if $G\not\simeq G[n,1,2]$, the proof of Theorem \ref{thm:1} is complete.

\begin{figure}[h!]
	\centering
	\includegraphics[width=3.75cm]{G2a.png}
	\hspace{1.5cm}
	\includegraphics[width=3.75cm]{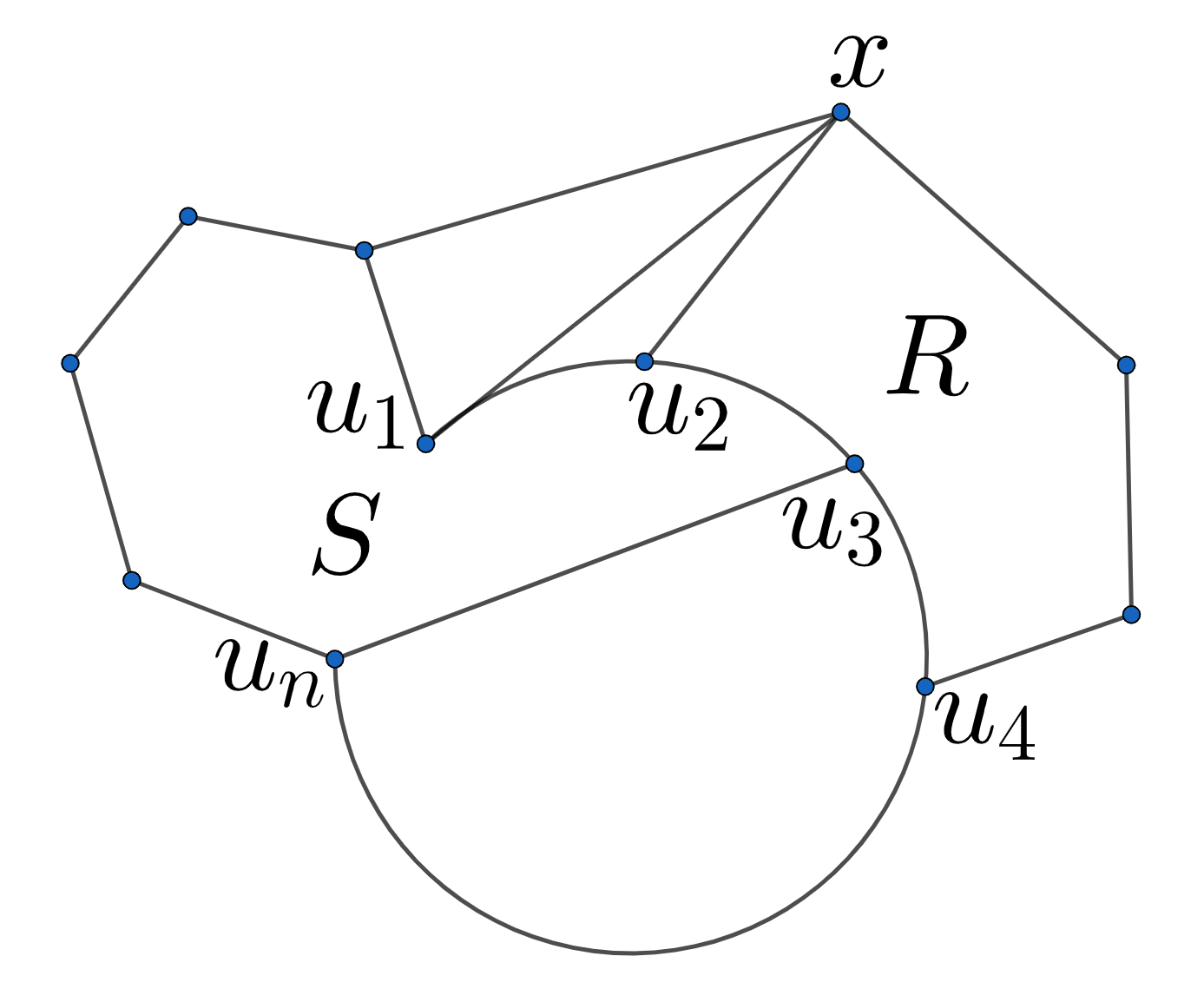}
	\caption{An example of $\U(n,1,3)$, transforming the subgraph of $G$ on the left to the subgraph of $G[n,1,3]$ on the right.}
	\label{fig:G3}
\end{figure}

Henceforth assume that \eqref{eq:cb} holds with $k=1$, $F_1$ is not a triangle, and $G\simeq G[n,1,2]$. We define
\begin{align*}
	&\U(n,1,3):G\to G[n,1,3],
	\\
	&G[n,1,3]=G+u_1x-xu_3+u_3u_n-u_nu_1,
\end{align*}
illustrated in Figure \ref{fig:G3}.

The transformed graph $G[n,1,3]$ is planar and of sequence $\sigma$. We point out that
\[G[n,1,2]\not\simeq G[n,1,3].\]
If $G[n,1,3]$ is $3$-connected, we have succeeded in proving Theorem \ref{thm:1}.

Assume henceforth that $G[n,1,3]$ is not $3$-connected. Similarly to the reasoning for the $3$-connectivity of $G[n,1,2]$, we deduce that $G[n,1,3]$ contains two distinct regions $R_1,R_2$ that share two non-adjacent vertices. We denote by $R$ the region of $G[n,1,3]$ obtained by adding $u_2$ to the vertices of $F_3$, and $S$ the region of $G[n,1,3]$ obtained by adding $u_2,u_3$ to the vertices of $F_n$, as in Figure \ref{fig:G2}. There are two cases to consider. First case, $R_1=R$ and $R_2$ contains $u_2$ and a vertex of $F_3$ distinct from $u_3,x$. The only regions $\neq R_1$ of $G[n,1,3]$ containing $u_2$ are $[u_1,u_2,x]$ and $S$, hence in this first case $R_2=S$. Second case, $R_1=S$ and $R_2$ contains one of $u_2,u_3$ and a vertex of $F_n$ distinct from $u_1,u_n$. We already know that $x\not\in F_n$, hence the only possibility is that $R_2=R$. Thereby, back to $G$, we have shown that in any case the faces $F_3$ and $F_n$ have a common vertex $z$, say, where $z\neq x$. In particular, $F_3$ is not a triangle.

To summarise: $G[n,1,2]$ is a polyhedron of same sequence as $G$. If it is not isomorphic to $G$, the proof is complete. If it is, we remark that $G[n,1,2]\not\simeq G[n,1,3]$, hence $G[n,1,3]$ is a simple, planar graph not isomorphic to $G$, and of same sequence $\sigma$. If $G[n,1,3]$ is $3$-connected, the proof is complete. Otherwise, in $G$ the faces $F_3$ and $F_n$ have a common vertex $z$, where $z\neq x$, and moreover $F_3$ is not a triangle.

Since $F_3$ is not a triangle, we can also define
\begin{align*}
	&\U(5,4,3):G\to G[5,4,3],
	\\
	&G[5,4,3]=G+u_4x-xu_3+u_3u_5-u_5u_4,
\end{align*}
and
\begin{align*}
	&\U(5,4,2):G\to G[5,4,2],
	\\
	&G[5,4,2]=G+u_4x-xu_2+u_2u_5-u_5u_4.
\end{align*}
We apply the above discussion of the proof of $G[n,1,2]$ being a polyhedron where the roles of $u_2,u_3$ are swapped, the roles of $u_1,u_4$ are swapped, and the roles of $u_n,u_5$ are swapped. Hence $G[5,4,3]$ is a polyhedron of same sequence as $G$. If it is not isomorphic to $G$, the proof is complete. If it is, we remark that $G[5,4,3]\not\simeq G[5,4,2]$, hence $G[5,4,2]$ is a simple, planar graph not isomorphic to $G$, and of same sequence $\sigma$.

We will now show that $G[5,4,2]$ is $3$-connected, completing the proof. Assume by contradiction that $G[5,4,2]$ is not 3-connected. We apply the above discussion, where we showed that if $G[n,1,3]$ is not 3-connected, then $F_3$ and $F_n$ have a common vertex $z\neq x$, but we swap the roles of $u_2,u_3$, of $u_1,u_4$, and of $u_n,u_5$. Hence in $G$ the faces $F_1$ and $F_4$ have a common vertex $z'$, where $z'\neq x$. Note that $F_3$ and $F_n$ have a common vertex $z$, and $F_1$ and $F_4$ have a common vertex $z'$, hence by planarity of $G$ we deduce that $z=z'$.

\begin{figure}[h!]
	\centering
	\includegraphics[width=4cm]{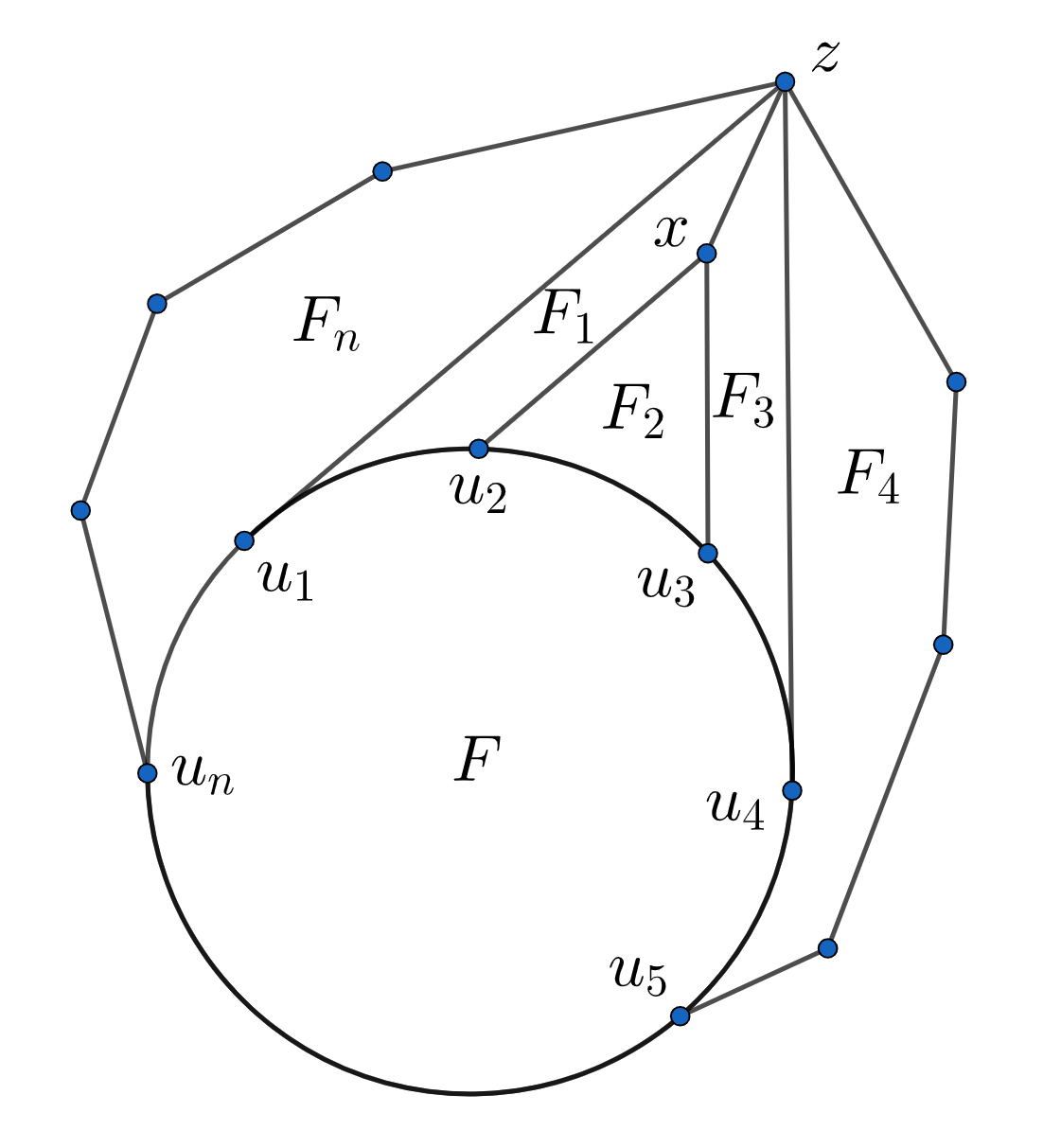}
	\hspace{1.5cm}
	\includegraphics[width=4cm]{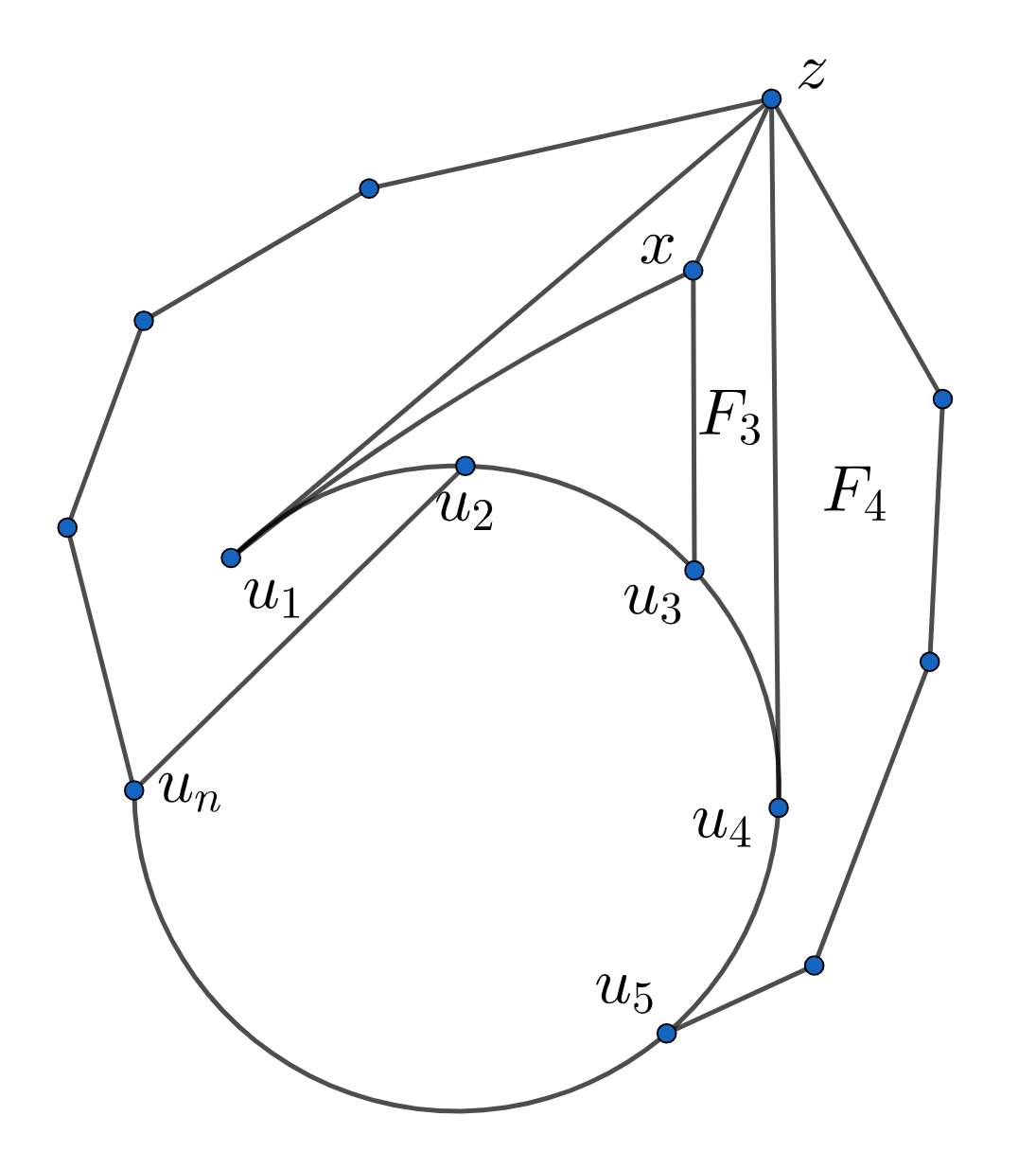}
	\caption{The situation when $F_1,F_3,F_4,F_n$ share a vertex $z$ (left), and the application of $\U(n,1,2)$ (right).}
	\label{fig:la}
\end{figure}

It follows that: $F_1,F_3$ share the vertices $x,z$, hence $xz\in E(G)$; $F_1,F_n$ share $u_1,z$, hence $u_1z\in E(G)$; and $F_3,F_4$ share $u_4,z$, hence $u_4z\in E(G)$. Therefore,
\[F_1=[u_1,u_2,x,z], \quad \text{ and } \quad F_3=[u_3,u_4,z,x].\]
In particular, all three vertices of the triangular face $F_2=[u_2,u_3,x]$ have degree $3$ in $G$. This situation is sketched in Figure \ref{fig:la}, left.

Finally, due to the isomorphism $\U(n,1,2)$, all vertices of the triangular face $[u_1,x,z]$ of $G[n,1,2]$ have degree $3$ in $G[n,1,2]$. Hence
\[\deg_G(z)=\deg_{G[n,1,2]}(z)=3.\]
However, the vertex $z$ cannot have degree $3$ in $G$, since the faces $F_4,F_n$ are distinct, contradiction. Hence, $G[5,4,2]$ is 3-connected, and thus a polyhedron.

Therefore, one of $G[n,1,2],G[n,1,3],G[5,4,3],G[5,4,2]$ is a polyhedron of same degree sequence as $G$. The proof of Theorem \ref{thm:1} is complete.

To summarise, given a non pyramidal polyhedron containing an $n$-gonal face for $n\geq 10$, we have constructed at least one non-isomorphic polyhedron of same degree sequence. This new graph is either one of $\varphi_i(G),\psi_i(G)$ for some $1\leq i\leq n$ (recall the definition \eqref{eqn:zeta} of these maps), or one of $G',G(2),G(3),G[n,1,2],G[n,1,3],G[5,4,3],G[5,4,2]$.

\paragraph{Acknowledgements.} The author wishes to thanks two anonymous referees for helpful comments and corrections on a previous version.
\\
R.~M.~was partially supported by Programme for Young Researchers `Rita Levi Montalcini' \textit{Discrete and Probabilistic Methods in Mathematics with Applications}, awarded to R.~M.

\bibliographystyle{abbrv}
\bibliography{bibgra}

\end{document}